\documentclass[12pt]{article}

\usepackage{amsmath}
\allowdisplaybreaks[4]

\usepackage[all]{xy}
\usepackage{amssymb}
\usepackage{amsthm}
\usepackage{hyperref}
\hypersetup{colorlinks=true,linkcolor=blue,citecolor=red}
\usepackage{amsmath}
\usepackage{amscd}
\usepackage{verbatim}
\usepackage{eurosym}
\usepackage{float}
\usepackage{color}
\usepackage{dcolumn}
\usepackage[mathscr]{eucal}
\usepackage[all]{xy}
\usepackage{hyperref}
\usepackage{mathrsfs}
\usepackage{amsmath}
\usepackage{amssymb}
\usepackage{amsfonts,ifpdf}
\usepackage{graphicx}
\usepackage{times}
\usepackage{float}
\usepackage{epstopdf}
\usepackage{cite}
\usepackage{youngtab}
\usepackage{ytableau}
\ytableausetup
{mathmode, boxsize=0.9em}

\newtheorem{pf}{proof}[section]
\theoremstyle{remark}

\makeatletter \@addtoreset{equation}{section} \makeatother
\makeindex 

\begin{document}

\begin{center}
{\Large\bf  Cubic Equations Through the Looking Glass of Sylvester
}
\end{center}

\vskip 6mm

\begin{center}
William Y.C. Chen\\[9pt]
Center for Applied Mathematics\\
Tianjin University\\
Tianjin 300072, P. R. China\\[9pt]
Email: chenyc@tju.edu.cn

\end{center}

\begin{center}
{\large\bf Abstract}
\end{center}

{\small One can hardly believe that there is still something to be said
about cubic equations. To dodge this doubt, we will instead try and
say something about Sylvester. He doubtless found a way of  solving
cubic equations. As mentioned by Rota, it was the
only method in this vein that  he could remember.
We realize that in the generic case Sylvester's magnificent approach aimed at reduced cubic
 equations boils down to an easy identity expressing a cubic polynomial
 as a sum of two third powers of linear forms. This leads to
 Cardano's formula  for cubic equations involving the
 third roots of unity.
}

\vskip 12mm

A special case of a remarkable discovery of
Sylvester \cite{Sylvester}
 states that in the generic
case, a cubic binary
form can be represented as a sum of
 two third powers of linear forms. This result
has been presented in the contexts of invariant theory of binary forms,
Waring's problem for binary forms,
the apolarity of polynomials and the umbral method, see \cite{ER, Kung, KR, Kung-Rota-Yan, Rota}.
However, it appears that no account of
this ingenious idea has
been given in full detail. Let us get to the point.

As is well known, by substituting $x$ with $x-a_1$,
a cubic polynomial
\begin{equation} \label{cp}
  x^3 +3a_1x^2 + 3a_2x +a_3
 \end{equation} can be written in the reduced
form
\begin{equation}
f(x)=x^3 - 3px+q.
\end{equation}
First, for the case $p=0$, the equation can be readily solved.
From now on, we assume that $p\not=0$.
For the moment, please do not ask why. Let us write $f(x)$ as
\begin{equation} \label{x-3-r-s}
 f(x) =x^3 - 3rsx + rs(r+s).
 \end{equation}
In   doing so, the parameters $r$ and $s$ are determined by the relations
 \[ rs =   p, \quad rs(r+s) =q.\]
 That is to say, $r$ and $s$ are the roots of the quadratic equation
 \begin{equation} \label{x-2-p-q}
 x^2 - {q\over p}x + p =0.
 \end{equation}

If $q^2 =4p^3$, the above equation (\ref{x-2-p-q})
 has  double roots, that is,
\[ r=s={q \over 2p}.\]
In this case, (\ref{x-3-r-s}) reduces to
\[f(x)=x^3 - 3r^2x + 2r^3,\]
and the following identity
\begin{equation}
x^3 - 3r^2x + 2r^3= (x-r)^2 (x+2r)
\end{equation}
serves as a key to the roots of $f(x)$.

Now we are  left with the generic case
when the quadratic equation
 (\ref{x-2-p-q}) has two distinct roots $r$ and $s$.
Under this circumstance, here emerges an identity:
\begin{equation}
 x^3 -3rs x +rs(r+s) =   {s \over s-r} (x-r)^3 - {r \over s-r } (x-s)^3,
 \end{equation}
which spells out why we choose to write $f(x)$ in the
form of (\ref{x-3-r-s}).
This relation also ensures that the cubic equation $f(x)=0$
can be solved with ease.

 The assumption that $p\neq 0$ implies that
  $s\neq 0$, and hence the cubic equation $f(x)=0$ can be reformulated as
\[(x-r)^3=\frac{r}{s}(x-s)^3.\]
Let $u_1,\,u_2,\,u_3$ be the three cubic roots of $r/s$.
Clearly, $u_1, u_2, u_3 \not=1$ since $r\not=s$.
Thus, the solutions of the equation $f(x)=0$
 are given by
\begin{equation}
x_j=\frac{r-su_j}{1-u_j},
\end{equation}
{where } $j=1,2,3.$

Consider the example in \cite{Bewersdorff}.
Let
\[f(x)=x^3-6x-6.\]
Then we are led to the relations $rs=2$ and $r+s=-3$.
Solve the quadratic equation
\[x^2+3x+2=0\]
to get
\[r=-2,\quad s=-1.\]
It follows that $f(x)$ can be expressed as
\begin{equation}\label{cubic}
f(x)=-(x+2)^3+2(x+1)^3,
\end{equation}
so that the cubic equation $f(x)=0$  takes the form
\[(x+2)^3=2(x+1)^3.\]
Let \[\omega=-\frac{1}{2}+\frac{\sqrt{3}}{2}i\]
be the  third root of unity, and let $u_1,\,u_2,\,u_3$ be the three cubic roots of $2$, namely,
for $j=0,1,2$,
\[u_j=\sqrt[3]{2}e^{j\frac{2\pi i}{3}}=\sqrt[3]{2}\omega^j.\]
Thus the solutions of the equation $f(x)=0$ are furnished by
\begin{eqnarray*}
 x_j=\frac{-2+u_j}{1-u_j}
 =\frac{\sqrt[3]{2}(\omega^j-\sqrt[3]{4})}{1-\omega^j\sqrt[3]{2}}
  =\omega^j\sqrt[3]{2}+\omega^{-j}\sqrt[3]{4},
\end{eqnarray*}
{where} $j=0,1,2$. It is no accident that the above solutions
 are in agreement with Cardano's formula
as presented in \cite{Bewersdorff}.

Last but not least, we should not forget that
we owe our thanks to Sylvester.

\vskip 6mm

\noindent
{\bf Acknowledgments.} The author is grateful to Joseph P.S. Kung, Peter Paule, Catherine H.F. Yan, Doron Zeilberger and the referees for their insightful comments and suggestions. This
work was supported by the National Science Foundation of China.


\begin{thebibliography}{99}

\bibitem{Bewersdorff}
Bewersdorff, J. (2006).
\emph{Galois Theory for Beginners: A Historical Perspective}. Providence, RI:
American Mathematical Society.

\bibitem{ER} Ehrenborg, R. and Rota, G.-C. (1993). Apolarity and
             canonical forms for homogeneous polynomials.
             European J. Combin. 14: 157--181.

\bibitem{KR} Kung, J.P.S. and Rota, G.-C. (1984). The theory of binary forms,
    Bull. Amer. Math. Soc. (N.S.) 10: 27--85.

\bibitem{Kung}  Kung, J.P.S. (1989). Canonical forms of binary forms: Variations
  on a theme of Sylvester, In: D. Stanton, Ed., \emph{Invariant Theory and
    Tableaux}, Springer-Verlag, pp. 46--58.

\bibitem{Kung-Rota-Yan} Kung, J.P.S., Rota, G.-C. and Yan, C.H.F. (2009).
\emph{Combinatorics: The Rota Way}. New York: Cambridge University Press.

\bibitem{Rota}  Rota, G.-C. (1999). Two turning points of invariant theory.
                  Math. Intelligencer 21: 20--27.

\bibitem{Sylvester} Sylvester, J.J. (1851). On a remarkable
discovery in the theory of canonical forms and of hyperdeterminants.
    Math. Mag. 2: 391--410.

\end{thebibliography}
\end{document}